\documentclass[10pt, reqno]{amsart}

\usepackage{amsfonts,amssymb,amsmath,amsthm,indentfirst}
\usepackage{fontenc,enumerate}

\usepackage{amsfonts}
\usepackage{amssymb}
\usepackage{latexsym}

\usepackage{tikz}
\usepackage{graphicx,graphics,epsfig}
\usepackage{epstopdf}

\usepackage[]{enumerate}
\usepackage{verbatim}
\usepackage{bbm}

\usepackage[]{enumerate}

\usepackage[sc]{mathpazo}          % Palatino with smallcaps as text font
\usepackage{eulervm}               % Euler math
\usepackage[scaled=0.86]{berasans} % Bera for san serif family
\usepackage[scaled=1]{inconsolata} % Inconsolata for fixed width
\usepackage[T1]{fontenc}

\providecommand{\keywords}[1]{\textbf{\textit{Keywords---}} #1}

\usepackage[final]{hyperref}
\hypersetup{colorlinks=true, linkcolor=blue, anchorcolor=blue, citecolor=red, filecolor=blue, menucolor=blue, pagecolor=blue, urlcolor=blue}

\newcommand{\norm}[1]{\left\Vert #1 \right\Vert}

\newcommand{\N}{\mathbb{N}}

\newcommand{\R}{\mathbb{R}}

\newtheorem{theorem}{Theorem}

\newtheorem{lemma}{Lemma}
\newtheorem{corollary}{Corollary}

\newtheorem{remark}{Remark}

\theoremstyle{definition}
\newtheorem{definition}{Definition}

\theoremstyle{remark}

\title[Well-posedness for MBE]{Well-posedness for a molecular beam epitaxy model}

\author [
L. Emerald] { 
Louis Emerald }

\author[ D. O. da Silva ]{Daniel Oliveira da Silva}

\author[A. Tesfahun]{Achenef Tesfahun}

\address{Department of Mathematics \\
Nazarbayev University \\
Qabanbai Batyr Avenue 53 \\
010000 Nur-Sultan \\
Republic of Kazakhstan}

\email{louisemerald76@gmail.com, achenef.tesfahun@nu.edu.kz}

\address{Department of Mathematics
\\
California State University, Los Angeles \\ 5151 State University Drive \\ Los Angeles, CA 90032
\\
United States}
 
   \email{ ddasilv4@calstatela.edu}

\keywords{Molecular beam epitaxy;  Well-posedness; Stochastic MBE}
\subjclass[2020]{35K25; 35K58; 35Q81}

\begin{document}

\begin{abstract} 
We study a general molecular beam epitaxy (MBE) equation modeling the epitaxial growth of thin films.  We show that, in the deterministic case, the associated Cauchy problem admits a unique smooth solution for all time, given initial data in the space $X_0=L^2(\R^d)\cap \dot W^{1,4}(\R^d)$ with $d=1, 2$.  This improves a recent result by Ag\'elas \cite{A15}, who established global existence in $H^3(\R^d)$.  Moreover, we investigate the local existence and uniqueness of solutions in the space $X_0$ for the stochastic MBE equation, with an additive noise that is white in time and regular in the space variable.
\end{abstract}

 \maketitle

\section{Introduction}
Molecular beam epitaxy (MBE) is a method of manufacturing very thin crystalline surfaces with well-defined orientations on substrates.  As the name suggests, this is accomplished by sending very precisely controlled beams of molecules into a vacuum system containing the substrate.  Originally discovered by Arthur and Lepore \cite{AL69}, MBE has become a leading method of manufacturing semiconductors used in various electronic components.  

To model the surfaces generated using MBE, various models of MBE have been proposed, such as:
\begin{itemize}
    \item the Edwards-Wilkinson model \cite{EW82}: $u_{t} - \alpha \Delta u = \eta$;
    \item the Kardar-Parisi-Zhang model \cite{KPZ86}: $u_{t} - \Delta u = \alpha | \nabla u |^{2} + \eta$;
    \item the Wolf-Villain model \cite{WV90}: $u_{t} + \alpha \Delta^{2} u = \eta$;
    \item the Zangwill model \cite{Z96}: $u_{t} + \alpha \Delta^{2} u = \beta \Delta ( |\nabla u|^{2} ) + \gamma \nabla \cdot \left( |\nabla u|^{2} \nabla u \right) + \eta$;
    \item the Rost-Krug model \cite{RK97}: $u_{t} + \Delta^{2}u + \nabla \cdot \left[ (1 - |\nabla u |^{2})\nabla u \right] = \eta$,
\end{itemize}
among others.  In these models, $u: \R^d \times \R_+ \rightarrow \R$ represents the deviation of the film height at the point $x$ from the mean film thickness at time $t$ and $\eta$ is a Gaussian white noise.  The constants $\alpha, \beta$, and $\gamma$ are generic constants which need not be the same in each model.

%An interesting phenomenon which has been observed experimentally is known as \emph{mound formation}.  Roughly speaking, this is the formation of pyramid-shaped structures on the surfaces generated using MBE. (For a scanning tunneling microscope image of these mounds, see \cite{KPM2000}.) To understand this phenomenon, various models of MBE have been proposed, such as,
%\begin{itemize}
%    \item the Edwards-Wilkinson model \cite{EW82}: $u_{t} - \alpha \Delta u = \eta$;
%    \item the Kardar-Parisi-Zhang model \cite{KPZ86}: $u_{t} - \Delta u = \alpha | \nabla u |^{2} + \eta$;
%    \item the Wolf-Villain model \cite{WV90}: $u_{t} + \alpha \Delta^{2} u = \eta$;
%    \item the Zangwill model \cite{Z96}: $u_{t} + \alpha \Delta^{2} u = \beta \Delta ( |\nabla u|^{2} ) + \gamma \nabla \cdot \left( |\nabla u|^{2} \nabla u \right) + \eta$;
%    \item the Rost-Krug model \cite{RK97}: $u_{t} + \Delta^{2}u + \nabla \cdot \left[ (1 - |\nabla u |^{2})\nabla u \right] = \eta$,
%\end{itemize}
%among others.  In these models, $u: \R^d \times \R_+ \rightarrow \R$ represents the deviation of the film height at the point $x$ from the mean film thickness at time $t$ and $\eta$ is a Gaussian white noise.  The constants $\alpha, \beta$, and $\gamma$ are generic constants which need not be the same in each model.

In an attempt to unify these models, Ag\'elas \cite{A15} introduced the equation
\begin{equation}\label{mbe}
u_{t} + \alpha \Delta^{2} u = N(u)  + \eta
%u_{t} + \Delta^{2} u = \alpha_1 \nabla\cdot\left(|\nabla u|^2\nabla u\right) - \alpha_2 \Delta\left(|\nabla u|^2\right) - \alpha_3 \Delta u + \alpha_4 |\nabla u|^2 + \eta
\end{equation}
where 
%$u: \R^d \times \R_+ \rightarrow \R$ and
\begin{equation}\label{N}
N(u)=\alpha_1 \nabla\cdot\left(|\nabla u|^2\nabla u\right) - \alpha_2 \Delta\left(|\nabla u|^2\right) - \alpha_3 \Delta u + \alpha_4 |\nabla u|^2 
\end{equation}
for given parameters $\alpha>0$ and $\alpha_j \ge 0$,  $j \in \{1, ... , 4 \}$.  By making the change of variable $t \mapsto t/\alpha$, we can from now on set $\alpha=1$.
For this generalized model \eqref{mbe}--\eqref{N}, Ag\'elas proved global well-posedness in $H^{s}$ for $s \geq 3$ under the assumption that $$\alpha_{1} > \alpha_{2}^{2}$$ in the deterministic case $\eta = 0$.  In the present work, we will extend Ag\'elas' result, reducing the required regularity to obtain global existence.  More specifically, it will be shown that the Cauchy problem for equation \eqref{mbe} in the deterministic case $\eta = 0$ is globally well-posed in the space $L^{2} \cap \dot{W}^{1, 4}$. Moreover, it is shown that the solution is smooth. These results will be stated more precisely in section \ref{introdet}, with the proofs contained in section \ref{wppf}.

An interesting fact about the literature on MBE is that the vast majority of results for any model of MBE are stated for the deterministic case $\eta = 0$, even though the model is inherently a stochastic one.  Some notable exceptions to this exist, namely the works of Bl\"omker et. al. \cite{B2015, B2001a, B2001b}.  To add to this small collection, we also present a local well-posedness result for the stochastic Ag\'elas model.  This will be stated precisely in section \ref{introstoch}, with the proof contained in section \ref{lwpstoch}.

Finally, we conclude with some comments regarding the phenomenon of \emph{coarsening}.  Ideally, surfaces generated by MBE should be perfectly flat.  However, the reality is that this is never the case.  In fact, it has been observed experimentally that exotic structures can appear on these surfaces; see \cite{KPM2000} for a scanning tunneling microscope image of these structures.  Thus, a major goal of MBE modeling is to quantify the ``roughness'' of these surfaces.  As is typical in engineering, deviations from the baseline are computed using the root-mean-square method.  Since solutions represent deviations from the average height, it is standard practice to define the \emph{coarseness} $C(t)$ of a solution at time $t$ to be its $L^{2}$ norm.  In section \ref{coarsening}, we make some comments on coarsening in these models, and point to some possible flaws in these models.

%A second goal of this work is to study the stochastic problem $\eta \neq 0$ for equation \eqref{mbe}.  Compared to the deterministic case, very few results exist for stochastic MBE equations; see, for example, the results contained in \cite{B2015, B2001a, B2001b}.  In the specific case of equation \eqref{mbe}, no results exist in the stochastic case.  In sections \ref{introstoch} and \ref{lwpstoch} below, the first result for equation \eqref{mbe} in the case $\eta \neq 0$ is presented.

%It should be noted that the global result in the deterministic case relies on a Gr\"onwall-type estimate for solutions to equation \eqref{mbe}.  This a priori estimate will have consequences regarding the \emph{coarseness} of solutions.  Coarsening is the phenomenon of imperfections appearing in the generated surfaces (such as the mounds mentioned previously).  In section \ref{coarsening}, we will discuss how coarseness is quantified, and the consequences of the a priori estimates related to coarsening.

%\vspace{2mm}
\subsection{Deterministic MBE equation ($\eta=0$)} \label{introdet}

In this section, we consider the Cauchy problem for the MBE equation
\begin{equation}\label{dmbe}
\left\{
\begin{aligned}
& u_{t} + \Delta^{2} u = N(u),
 \\
& u(x,0) = u_0(x),
\end{aligned}
\right.
\end{equation}
%\vspace{2mm}
where $N(u)$ is as in equation \eqref{N}.  To state our result, we introduce the notation
$$
X_0=L^2(\mathbb{R}^d) \cap \dot W^{1,4}(\mathbb{R}^d).
$$
%\vspace{2mm}
Our main result concerning \eqref{dmbe} is:
\begin{theorem}[Unconditional Well-posedness]\label{thm-wp}
Assume that $u_{0} \in X_{0}$.
\begin{enumerate}[(i)]
\item \label{wpi}   Let $d=1,2, 3$. Then  there exists $T>0$ and a unique solution $u \in C\left([0, T); X_0\right)$ of \eqref{dmbe}. Moreover, the solution is smooth, i.e., 
 \begin{equation}\label{smooth-loc}
 u \in C^\infty\left(  (0,T); C^\infty(\mathbb{R}^d) \right).
 \end{equation}
 
\item  \label{wpii} Let $d=1,2$ and $\alpha_1 > \alpha_2^2$. Then there exists a unique global smooth solution $u \in C^{\infty}((0,+\infty);C^{\infty}(\mathbb{R}^d))$ of \eqref{dmbe}.
 \end{enumerate}
 
\end{theorem}
\begin{remark}
The local result can be shown to hold in any dimension, but for $d\ge 4$, it is necessary to modify the initial data space to $$X_0 = L^2(\mathbb{R}^d) \cap \dot W^{1,p}(\mathbb{R}^d)$$ with $p>d.$
\end{remark}
The local well-posedness result  in Theorem \ref{thm-wp}\eqref{wpi} is proved by applying a fixed point argument to the mild formulation of 
 \eqref{dmbe}, 
\begin{equation}\label{dmbe-I}
	u(t)=  S(t)  u_0
	+  \int_0^t S(t-s)  N(u (s) )\; ds
\end{equation}
in the space
$X_T=  C \left( [0,T);  X_0\right)$
 with associated norm
$
\norm{u}_{X_T} = \norm{u}_{L_T^\infty L_x^2} + \norm{ \nabla u}_{L_T^\infty  L_x^4},
$
where we use the notation $L^q_T Y: = L^q\left( [0, T); Y\right)$ for a Banach space $Y$.  
Here $S(t)= e^{-t \Delta^2}$ is the linear propagator with Fourier symbol $e^{-t |\xi|^4}$, i.e., 
$$
\mathcal S(t) u_0 = \mathcal F^{-1} \left[  e^{-t |\xi|^4} \widehat{f}\right]= (K_t \ast u_0) (x),
$$
where
$$
K_t (x) = \int_{\R^d} e^{-t |\xi|^4} e^{ix\xi} \, d\xi$$ is the kernel.  As is well-known, it satisfies the time-decay estimate 
\begin{equation}\label{Kdecay-est}
   \norm{ \partial_x^{\alpha}  K_t}_{L_x^p (\R^d)} \le C t^{-\frac d 4\left(1-\frac 1p \right) -\frac {|\alpha|}  4}
   \end{equation}
    for all $\alpha \in \mathbb{N}^d$ , $p\ge 1$ and some constant $C=C(\alpha, p)>0$. For the proof, see e.g., \cite[Lemma 2.1]{LWZ2007}. 
    Similarly, one can show that \begin{equation}\label{Kdecay-est-s}
   \norm{ |D|^s K_t}_{L_x^p (\R^d)} \le C t^{-\frac d 4\left(1-\frac 1p \right) -\frac {s}  4} 
   \end{equation}
for all $s\ge 0$. 

By combining \eqref{Kdecay-est} and \eqref{Kdecay-est-s} with Young's inequality, we can obtain the following estimates for $1\le
   r \leq p\le \infty$:
\begin{equation}\label{decay-est integer derivatives}
   \norm{ \partial_x^{\alpha} S(t) u_0}_{L_x^p (\R^d)} \le C t^{-\frac d 4\left(\frac1r-\frac 1p \right) -\frac{|\alpha|}{4}}   \norm{  u_0}_{L_x^r (\R^d)} 
   \end{equation}
   and
      \begin{equation}\label{decay-est real derivatives}
   \norm{ |D|^{s} S(t) u_0}_{L_x^p (\R^d)} \le C t^{-\frac d 4\left(\frac1r-\frac 1p \right) -\frac{s}{4}}   \norm{  u_0}_{L_x^r (\R^d)}.
   \end{equation}

The proof of the local existence result relies on estimate \eqref{decay-est integer derivatives}. 
To obtain the smoothness result \eqref{smooth-loc}  in the spatial variable $x$, we use \eqref{decay-est integer derivatives}--\eqref{decay-est real derivatives} as well as 
 the standard product estimate for the Sobolev spaces (see e.g., Corollary 3.16 of \cite{Tao2001}),
    \begin{equation}\label{prodest}
            \norm{fg}_{H^s} \lesssim \norm{f}_{H^{s_1}} \norm{g}_{H^{s_2}} \quad \text{if} \quad \begin{cases} 
        &s_1 + s_2 \ge 0, 
        \\
        &s\le s_1,  \ s\le s_2
        \\
       & s < s_1 + s_2 - \frac{d}{2}
        \end{cases}
    \end{equation}
    for all $f \in H^{s_1}(\R^d), g \in H^{s_2}(\R^d)$.
Smoothness in the time variable $t$ then follows from the PDE \eqref{dmbe} by a standard bootstrap argument. The global existence result in  Theorem \ref{thm-wp}\eqref{wpii} will follow from 
\eqref{smooth-loc} and the global existence result in $H^3$ due to Ag\'elas \cite{A15}.

%\vspace{2mm}

%\vspace{2mm}
\subsection{Stochastic MBE equation ($\eta\neq 0$)}\label{introstoch}
%\vspace{2mm}

In addition to the deterministic case, we also consider the Cauchy problem for the stochastic MBE equation
\begin{equation}\label{smbe}
\left\{
\begin{aligned}
& u_{t} + \Delta^{2} u = N(u) + \eta
 \\
& u(x,0) = u_0(x),
\end{aligned}
\right.
\end{equation}
  where
$
\eta(x,t)
$
is a real-valued Gaussian process with correlation function
$$
\mathbb E \left[ \eta(x,t)  \eta(y, s)  \right]= c(x,y) \delta (t-s)
$$
for some function $c$.
The case $c(x,y) =\delta(x-y)$ corresponds to a space-time white noise. However, this is difficult to treat mathematically on unbounded domains, and therefore smoother correlation function will be considered here.

 In order to give a mathematical definition of $\eta$ and write \eqref{smbe} as an It\^{o} integral, we
introduce a probability space $( \Omega, \mathcal F, \mathbb P)$ endowed with a filtration $\{ \mathcal F_t \}_{t\ge 0}$, and a sequence $\{ \beta_k \}_{k\in \N}$ of independent real valued Brownian motions on $\R_+$ associated to the filtration  $\{ \mathcal F_t \}_{t\ge 0}$.
Given an orthonormal basis $\{e_k\}_{k\in \N}$ of $L^2(\R^d)$ and a symmetric 
linear operator $\phi$ on  $L^2(\R^d)$, the process
$$
W_\phi(x,t, \omega)=\sum_{k\in \N} \phi e_k (x)  \cdot  \beta_k(t, \omega) \quad,  x\in \R, \ \  t\ge 0,  \ \ \omega \in \Omega
$$
is a Wiener process on $L^2(\R^d)$ with covariance operator $\phi \phi^\ast$.

Note that $\phi=I$ corresponds to a space-time white noise with $c(x,y) =\delta(x-y)$. If $\phi$ is defined through a kernel $L$, which means that for any square integrable function $f$,
$$
Lf(x)= \int_{\R^d} L(x-y) f(y)\, dy, \qquad  f\in L^2(\R^d),
$$
then 
$$
c(x,y)= \int_{\R^d}  L(x, z) L(y, z)\, dz.
$$
In general, we assume that 
 $\phi$ is a Hilbert-Schmidt operator from $L^2(\R^d)$ to $H^s(\R^d)$ for appropriate values of $s$, endowed with the norm
$$
\norm{\phi}^2_{HS(L^2: H^s)}= \sum_{k\in \N} \norm {\phi e_k}_{H^s}^2.
$$

Now,  we set 
$$
\eta= \frac{dW_\phi}{dt}.
$$
Then the  equivalent It\^{o}  equation to \eqref{smbe} becomes
\begin{equation}\label{smbe-ito}
\left\{
\begin{aligned}
&du =- \Delta^{2} u  dt +N(u) dt  + dW_\phi
\\
& u(x,0) = u_0(x).
\end{aligned}
\right.
\end{equation}
  We shall construct a solution $u$ satisfying the following mild
formulation of \eqref{smbe-ito}:
\begin{equation}\label{ismbe}
	u(t)=   S(t) u_0
	+ \int_0^t S(t-s) N(u(s))\; ds + \int_0^t S(t-s)  dW_\phi(s) .
\end{equation}
The last term on the right-hand side of \eqref{ismbe} represents the effect of the stochastic forcing and is called the \emph{stochastic convolution}, which we denote by
$Z$, i.e.,
\begin{equation}\label{stconv}
\begin{split}
 Z(t) &=\int_0^t S(t-s)  dW_\phi(s) 
 \\
 &=\sum_{k\in \N} \int_0^t S(t-s) (\phi e_k)(x)  \; d \beta_k(s) .
 \end{split}
\end{equation}

\begin{comment}
This can be defined path-wise without Ito-calculus using the integration by parts formula:
\begin{align*}
 Z(t) &= W_\phi(t)-  \Delta^2 \int_0^t S(t-s)    W_\phi(s) \, ds .
\end{align*}
So formally $Z$ solves the IVP
$$
\partial_t Z + \Delta^2 Z = \partial_t W, \quad Z(0)=0.
$$
\end{comment}

From the mild formulation \eqref{ismbe}-\eqref{stconv} we see that any solution $u$ can be written as $u=v+Z$.
We then study a fixed point problem for the residual term 
 $v=u-Z$:
\begin{equation}\label{smbe-I-resd}
v(t)=   S(t)  u_0
	+  \int_0^t S(t-s) N(v+Z)  (s)\; ds.
\end{equation}

\begin{theorem}\label{thm-slwp}
Let  $1\le d\le 3$,  $u_0 \in  X_0$ and 
and $\phi \in HS(L^2; L^2)$. Then there exists a stopping time $T=T(\omega)>0$ 
and a unique solution $v \in C\left([0, T]; X_0\right)$ of \eqref{smbe-I-resd} almost surely for $\omega \in \Omega$.
\end{theorem}

%\vspace{2mm}

\begin{lemma}[A priori estimate] \label{lm-ap} Let $\alpha_1 = \alpha_3 = 1$, $\alpha_2 = \alpha_4 = 0$.
 If $v$ is a smooth solution of \eqref{smbe-I-resd} on $[0, T]$, then
\begin{align*}
 \norm{v(t)}_{L_x^2}^2 + 2 \int_0^t \left[ \norm{ \Delta v (s)}_{L_x^2}^2 +  \norm{ \nabla v (s)}_{L_x^4}^4\right ] \, ds \le 2\left[ \norm{u_0}^2_{L_x^2}  +  C \right] e^{T}  
\end{align*}
almost surely for $\omega \in \Omega$, where $C=C \left (T, \norm{\phi}_{HS(L^2: L^2)}\right)>0$.

\end{lemma} 
\begin{remark}
    It is also possible to derive analogous a priori estimates for the general case $\alpha_j\ge 0$ ($j=1,2, 3,4$).  However, we do not pursue this here.
\end{remark}

\section{Proof of Theorem \ref{thm-wp} }\label{wppf}

\subsection{Local existence of solutions}

Define the mapping
$\Gamma=\Gamma_{u_0}$ by
\begin{align}\label{mbe-fI}
	\Gamma u(t)=  S(t)  u_0
	+  \int_0^t S(t-s)  N(u(s))\; ds.
\end{align}
Given initial data with norm
 $\norm{u_0}_{X_0} \le R$, we show that $\Gamma$ is a contraction in the ball
$$
B_{ R, T}= \left\{  u \in X_T : \  \norm{u}_{X_T}  \le 2 R \right\},
$$
for sufficiently small $T=T(R)>0$.  For this, it suffices to prove 
\begin{align}
\label{Gumap}
\norm{\Gamma u}_{X_T} & \le   2R 
\\
\label{Gumap-d}
\norm{\Gamma u- \Gamma v}_{X_T}  < &\norm{u- v}_{X_T}
\end{align}
 for $u, v\in  B_{ R, T}$. Therefore, the map $\Gamma$ has a unique fixed point $u\in _{ R, T}$
 solving the integral equation \eqref{dmbe-I}.  Continuous dependence on the initial data can be shown in a similar way.

As no additional difficulty arises when estimating the lower order terms, for the sake of clarity we set $\alpha_3 = 0$ and $\alpha_4 = 0$. Moreover, as the value of the parameters $\alpha_1$ and $\alpha_2$ are not of significant importance, we set $\alpha_1 = \alpha_2 = 1$.  
So we take $N(u)=\nabla\cdot\left(|\nabla u|^2\nabla u\right) - \Delta\left(|\nabla u|^2\right)$.  Then by \eqref{decay-est integer derivatives},
\begin{align*}
	 \norm{ \Gamma u(t)}_{L_x^2} &\le  \norm{ S(t)  u_0 }_{L_x^2} + \int_0^t  \norm{ \nabla S(t-s) \cdot |\nabla u(s)|^2 \nabla u (s)}_{L_x^{2}}   \; ds \\
     & \quad + \int_0^t \norm{ \Delta S(t-s) |\nabla u(s)|^2 }_{L_x^{2}} \; ds 
	\\
	 & \le \norm{ u_0}_{L_x^2} +  \int_0^t (t-s)^{-\frac{d+4}{16}}  \norm{|\nabla u(s)|^3}_{L_x^{\frac{4}{3}}} \; ds \\
      & \quad + \int_0^t (t-s)^{-\frac{1}{2}}  \norm{|\nabla u(s)|^2}_{L_x^{2}}   \; ds 
        \\
	  & \le \norm{u_0}_{L_x^2} + C \left[ T^{\frac{12-d}{16}} \norm{\nabla u}^3_{ L^\infty_T L_x^4} + T^{\frac{1}{2}} \norm{\nabla u}^2_{ L^\infty_T L_x^4}  \right].
\end{align*} 

Similarly, for $\alpha \in \mathbb{N}^d$, such that $|\alpha| = 1 $, we have
\begin{align*}
	 \norm{ \partial_x^{\alpha} \Gamma u(t)}_{L_x^4} &\le  \norm{  S(t) \partial_x^{\alpha} u_0 }_{L_x^4}   + \int_0^t  \norm{ \nabla \partial_x^{\alpha} S(t-s) \cdot |\nabla u(s)|^2 \nabla u (s)}_{L_x^4}   \; ds 
      \\
     &\hspace{0.3cm} + \int_0^t   \norm{ \Delta \partial_x^{\alpha} S(t-s) |\nabla u(s)|^2 }_{L_x^{4}} \; ds 
	\\
	 & \le \norm{ u_0}_{\dot W^{1,4}} +  \int_0^t (t-s)^{-\frac{d}{2p} - \frac{1}{2} }  \norm{|\nabla u(s)|^3}_{L_x^{\frac{4}{3}}} \; ds \\
      & \quad + \int_0^t (t-s)^{-\frac{d}{16} - \frac{3}{4}}  \norm{|\nabla u(s)|^2}_{L_x^{\frac{p}{2}}}   \; ds 
	 \\
	  & \le \norm{ u_0}_{\dot W^{1,4}} + C  \left[ T^{\frac{1}{2}-\frac{d}{8}} \norm{u}^3_{ L^\infty_T \dot W^{1,4}} + T^{\frac{1}{4}-\frac{d}{16}} \norm{u}^2_{ L^\infty_T \dot W^{1,4}}  \right].
\end{align*}
From the above estimates, we conclude that 
\begin{equation*}
\begin{split}
\norm{\Gamma u}_{X_T} & \le \norm{u_0}_{ X_0}   + C  T^{\frac{1}{2}-\frac{d}{8}}  \left[ \norm{u}^2_{X_T} + \norm{u}^3_{X_T}\right]
\\
&\le 2R,
\end{split}
\end{equation*}
provided that
$$
T\le \left[ 4CR(1+2R)\right]^{-\frac8{4-d}}.
$$

In a similar way, we can derive the difference estimate 
\begin{equation*}
\begin{split}
\norm{\Gamma u-\Gamma v}_{X_T}  & \le      C T^{\frac{1}{2}-\frac{d}{8}}  \left[ \norm{u}_{X_T} +  \norm{v}_{X_T} + \norm{u}^2_{X_T} +  \norm{v}_{X_T}^2\right] \norm{ u- v}_{X_T} 
\\
&< \norm{ u- v}_{X_T},
\end{split}
\end{equation*}
which holds when
\[
T < \left[C\left( 4R + 8R^{2} \right)\right]^{-\frac{8}{4 - d}}.
\]
It follows that $\Gamma u$ is a contraction on $X_{T}$, from which existence and uniqueness of solutions follows.

\subsection{Regularity of the local solution}

Our next goal is to prove the regularity result in part (i) of Theorem \ref{thm-wp}.  As a first step, we prove the following lemma:
\begin{lemma}\label{Lemma regularity}
    Let $u \in X_T$ be a solution of \eqref{dmbe-I} associated with the initial condition $u_0 \in X_0$. Then $u \in L^{\infty}((0,T];C^{\infty}(R^d))$.
\end{lemma}
\begin{proof}

For the sake of clarity, like in the proof of Theorem \ref{thm-wp}, we set $\alpha_1 = \alpha_2 = 1$ and $\alpha_3 = \alpha_4 = 0$.
    Using the semi-group property of the linear propagator $S$, we have that for any $0 \leq \overline{t} < t < T$, 
    \begin{equation}\label{mild formulation with time shift}
    \begin{split}
        u(t) = S(t-\bar{t})u(\bar{t}) &+ \int_{\bar{t}}^t \nabla \cdot S(t-s) |\nabla u(s)|^2 \nabla u(s) \mathrm{d}s 
        \\
        &+
        \int_{\bar{t}}^t \Delta S(t-s) |\nabla u(s)|^2 \mathrm{d}s.
        \end{split}
    \end{equation}
    The proof consists of four steps. \\
   
   \noindent {\bf Step 1.}  Taking first $\bar{t} = 0$, applying the operator $|D|^{\alpha}$ to this equation \eqref{mild formulation with time shift} for $\alpha \in \mathbb{R}$, then taking the $L^2$ norm and using Hausdorff's inequality, we get
    \begin{multline*}
        \norm{|D|^{\alpha} u(t)}_{L^2} \leq \norm{|D|^{\alpha}S(t)u_0}_{L^2} +  \int_{0}^t \norm{|D|^{\alpha}\nabla S(t-s) \cdot |\nabla u(s)|^2 \nabla u(s)}_{L^2} \mathrm{d}s
        \\
        + \int_{0}^t \norm{|D|^{\alpha} \Delta S(t-s) |\nabla u(s)|^2}_{L^2} \mathrm{d}s. 
    \end{multline*}
    Using the decay estimates \eqref{decay-est real derivatives} of $S$, we get
    \begin{align*}
        \norm{|D|^{\alpha} u(t)}_{L^2} &\lesssim t^{-\frac{\alpha}{4}} \norm{u_0}_{L^2} + \int_{0}^t (t-s)^{-\frac{d + 4\alpha + 4}{16}} \norm{|\nabla u(s)|^2 \nabla u(s)}_{L^{\frac{4}{3}}} \mathrm{d}s 
        \\ 
        &\hspace{0.3cm} + \int_{0}^t (t-s)^{-\frac{\alpha+2}{4}} \norm{|\nabla u(s)|^2}_{L^2} \mathrm{d}s
        \\
        &\lesssim t^{-\frac{\alpha}{4}} \norm{u_0}_{L^2} + \int_{0}^t (t-s)^{-\frac{d+4\alpha+4}{16}} \norm{\nabla u(s)}_{L^4}^3 \mathrm{d}s 
        \\ 
        &\hspace{0.3cm} + \int_{0}^t (t-s)^{-\frac{\alpha+2}{4}} \norm{\nabla u(s)}_{L^4}^2 \mathrm{d}s.
    \end{align*}
    So that, if $\alpha < 2$, we get
    \begin{align*}
        \norm{|D|^{\alpha} u(t)}_{L^2} \lesssim t^{-\frac{\alpha}{4}} \norm{u_0}_{X_0} + t^{\frac{2 -\alpha}{4}} \norm{u}_{X_T}^2 + t^{\frac{12 - d - 4\alpha}{16}}\norm{u}_{X_T}^3.
    \end{align*}
So for now, we have $u \in L^{\infty}_{(0,T]} H^{2-}(\mathbb{R}^d)$, since $u\in X_T$ (by Theorem \ref{thm-wp}). \\

    \noindent {\bf Step 2.} This step is required only in the case $d=2$. In the case $d=1$, one can go directly to the Step 3 of the proof. Let us apply the operator $|D|^{\alpha + \frac{1}{2}}$, with $\alpha < 2$, to \eqref{mild formulation with time shift}.
    \begin{align*}
        \norm{|D|^{\alpha + \frac{1}{2}} u(t)}_{L^2} \leq &\norm{|D|^{\alpha + \frac{1}{2}}S(t) u_0}_{L^2} \\
        & + \int_{0}^t \norm{|D|^{\alpha+\frac{1}{2}}\nabla S(t-s) \cdot \big(|\nabla u(s)|^2 \nabla u(s)\big)}_{L^2} \mathrm{d}s
        \\
        & + \int_{0}^t \norm{|D|^{\alpha} \Delta S(t-s) |D|^{\frac{1}{2}}\big(|\nabla u(s)|^2\big)}_{L^2} \mathrm{d}s.
    \end{align*}
    Using again the decay estimates \eqref{decay-est real derivatives}, we have 
    \begin{align*}
        \norm{|D|^{\alpha + \frac{1}{2}} u(t)}_{L^2} &\lesssim t^{-\frac{\alpha+\frac{1}{2}}{4}} \norm{u_0}_{L^2} + \int_{0}^t (t-s)^{-\frac{2 + \alpha}{4}} \norm{|\nabla u(s)|^2 \nabla u(s)}_{L^{\frac{4}{3}}} \mathrm{d}s 
        \\ 
        &\hspace{0.3cm} + \int_{0}^t (t-s)^{-\frac{\alpha + 2}{4}} \norm{|\nabla u(s)|^2}_{H^{\frac{1}{2}}} \mathrm{d}s.
    \end{align*}
    Now, in order to estimate the last term of the right hand side, we use \eqref{prodest}, with $s = \frac{1}{2}, s_1 = s_2 = \frac{7}{8}$. Then $s \leq s_j$ and $s_1 + s_2 - 1 = \frac{3}{4} > s$. We get
    \begin{align*}
        \norm{|D|^{\alpha + \frac{1}{2}} u(t)}_{L^2} &\lesssim t^{-\frac{2\alpha+1}{8}} \norm{u_0}_{X_0} + t^{\frac{2 - \alpha}{4}} \norm{u}_{X_T}^3 
        \\ 
        &\hspace{0.3cm} + \int_{0}^t (t-s)^{-\frac{\alpha + 2}{4}} \norm{u(s)}_{H^{1+\frac{7}{8}}}^2 \mathrm{d}s
        \\
        &\lesssim t^{-\frac{2\alpha+1}{8}}\norm{u_0}_{X_0} +  t^{\frac{2 - \alpha}{4}}\left[ \norm{u}_{L^{\infty}_{(0,T]}H^{\alpha}}^2 + \norm{u}_{X_T}^3\right],
    \end{align*}
    where in the last line we choose $\alpha$ sufficiently large so that $2 > \alpha > 1 + \frac{7}{8}$.
    Now, we know that $u \in L^{\infty}_{(0,T]} H^{2 + \frac{1}{2}-}(\mathbb{R}^d)$. \\

\noindent {\bf Step 3.} Let $s_0 > \frac{d}{2}$ and apply $|D|^{\alpha + s_0}$ to \eqref{mild formulation with time shift}.  We then have
    \begin{align*}
        \norm{|D|^{\alpha + s_0} u(t)}_{L^2} \leq &\norm{|D|^{\alpha}S(t- \bar{t})|D|^{s_0}u(\bar{t})}_{L^2} \\
        & +  \int_{\bar{t}}^t \norm{|D|^{\alpha}\nabla S(t-s) \cdot |D|^{s_0}\big(|\nabla u(s)|^2 \nabla u(s)\big)}_{L^2} \mathrm{d}s
        \\
        & + \int_{\bar{t}}^t \norm{|D|^{\alpha} \Delta S(t-s) |D|^{s_0}\big(|\nabla u(s)|^2\big)}_{L^2} \mathrm{d}s,
    \end{align*}
    where $0 < \bar{t} \leq t < T$.
    Then, using the decay estimates \eqref{decay-est real derivatives}, we get
    \begin{align*}
        \norm{|D|^{\alpha + s_0} u(t)}_{L^2} &\lesssim (t-\bar{t})^{-\frac{\alpha}{4}} \norm{|D|^{s_0} u(\bar{t})}_{L^2} + \int_{\bar{t}}^t (t-s)^{-\frac{\alpha+1}{4}} \norm{|\nabla u(s)|^2 \nabla u(s)}_{H^{s_0}} \mathrm{d}s 
        \\ 
        &\hspace{0.3cm} + \int_{\bar{t}}^t (t-s)^{-\frac{\alpha + 2}{4}} \norm{|\nabla u(s)|^2}_{H^{s_0}} \mathrm{d}s
        \\
        &\lesssim (t-\bar{t})^{-\frac{\alpha}{4}} \norm{u}_{L^{\infty}_{(0,T]}H^{s_0+1}} + (t-\bar{t})^{\frac{2-\alpha}{4}} \norm{u}_{L^{\infty}_{(0,T]} H^{s_0+1}}^2 + (t-\bar{t})^{\frac{3-\alpha}{4}}\norm{u}_{L^{\infty}_{(0,T]}H^{s_0+1}}^3 , 
    \end{align*}
    where we used the algebra property of $H^{s_0}$, and where we take $s_0$ sufficiently small so that $s_0 + 1 = \frac{d}{2} + 1 + \varepsilon < 2$ when $d = 1$, and $s_0 + 1 = \frac{d}{2} + 1 + \varepsilon < 2 + \frac{1}{2}$ when $d = 2$. Now, we know that $u \in L^{\infty}_{(0,T]} H^{2 + \frac{d}{2}}(\mathbb{R}^d)$. \\
    
\noindent {\bf Step 4.} Let $n \geq d$, let us suppose that $u \in L^{\infty}_{(0,T]} H^{2 + \frac{n}{2}}(\mathbb{R}^d)$, which holds when $n = d$.  Let us prove that $u \in L^{\infty}_{(0,T]} H^{2 + \frac{n + 1}{2}}(\mathbb{R}^d)$.  From \eqref{mild formulation with time shift}, we can apply the operator $|D|^{\gamma_1 + \gamma_2}$ where $\gamma_1 = \frac{3}{2}$ and $\gamma_2 = 2 + \frac{n}{2} - 1$ (so that $\gamma_1 < 2$ and $\gamma_2 + 1 = 2 + \frac{n}{2}$ ) to obtain
    \begin{align*}
        \norm{|D|^{\gamma_1 + \gamma_2} u(t)}_{L^2} \leq &\norm{|D|^{\gamma_1}S(t- \bar{t})|D|^{\gamma_2}u(\bar{t})}_{L^2} \\
        & +  \int_{\bar{t}}^t \norm{|D|^{\gamma_1}\nabla S(t-s) \cdot |D|^{\gamma_2}\big(|\nabla u(s)|^2 \nabla u(s)\big)}_{L^2} \mathrm{d}s
        \\
        & + \int_{\bar{t}}^t \norm{|D|^{\gamma_1} \Delta S(t-s) |D|^{\gamma_2}\big(|\nabla u(s)|^2\big)}_{L^2} \mathrm{d}s.
    \end{align*}
    Then, using the decay estimates \eqref{decay-est real derivatives}, we get
    \begin{align*}
        \norm{|D|^{\gamma_1 + \gamma_2} u(t)}_{L^2} \lesssim & (t-\bar{t})^{-\frac{\gamma_1}{4}} \norm{|D|^{\gamma_2} u(\bar{t})}_{L^2} \\
        & + \int_{\bar{t}}^t (t-s)^{-\frac{\gamma_1+1}{4}} \norm{|\nabla u(s)|^2 \nabla u(s)}_{H^{\gamma_2}} \mathrm{d}s 
        \\ 
        &\hspace{0.3cm} + \int_{\bar{t}}^t (t-s)^{-\frac{\gamma_1 + 2}{4}} \norm{|\nabla u(s)|^2}_{H^{\gamma_2}} \mathrm{d}s 
        \\
        &\lesssim (t-\bar{t})^{-\frac{\gamma_1}{4}} \norm{u}_{L^{\infty}_{(0,T]}H^{\gamma_2+1}} \\
        & + (t-\bar{t})^{\frac{2 - \gamma_1}{4}} \norm{u}_{L^{\infty}_{(0,T]} H^{\gamma_2+1}}^2 +  (t-\bar{t})^{\frac{3 - \gamma_1}{4}}\norm{u}_{L^{\infty}_{(0,T]}H^{\gamma_2+1}}^3 , 
    \end{align*}
    where we used the algebra property of $H^{\gamma_2}$.  This completes the proof.
    \end{proof}

\begin{remark}
    The proof of Lemma \ref{Lemma regularity} can be adapted to the case $d = 3$. 
 Step 1 would be the same.  Step 2 is modified in the following way:

    \noindent {\bf Step 2'.} Let us apply the operator $|D|^{\alpha + \frac{1}{4}}$, with $\alpha < 2$, to \eqref{mild formulation with time shift}.
    \begin{align*}
        \norm{|D|^{\alpha + \frac{1}{4}} u(t)}_{L^2} \leq &\norm{|D|^{\alpha + \frac{1}{4}}S(t) u_0}_{L^2} +  \int_{0}^t \norm{|D|^{\alpha+\frac{1}{4}}\nabla S(t-s) \cdot \big(|\nabla u(s)|^2 \nabla u(s)\big)}_{L^2} \mathrm{d}s
        \\
        + &\int_{0}^t \norm{|D|^{\alpha} \Delta S(t-s) |D|^{\frac{1}{4}}\big(|\nabla u(s)|^2\big)}_{L^2} \mathrm{d}s.
    \end{align*}
    Using the decay estimates \eqref{decay-est real derivatives}, we have 
    \begin{align*}
        \norm{|D|^{\alpha + \frac{1}{4}} u(t)}_{L^2} &\lesssim t^{-\frac{\alpha+\frac{1}{4}}{4}} \norm{u_0}_{L^2} + \int_{0}^t (t-s)^{-\frac{2 + \alpha}{4}} \norm{|\nabla u(s)|^2 \nabla u(s)}_{L^{\frac{4}{3}}} \mathrm{d}s 
        \\ 
        &\hspace{0.3cm} + \int_{0}^t (t-s)^{-\frac{\alpha + 2}{4}} \norm{|\nabla u(s)|^2}_{H^{\frac{1}{4}}} \mathrm{d}s.
    \end{align*}
    Now, let us use \eqref{prodest} for the term $I := \int_{0}^t (t-s)^{-\frac{\alpha + 2}{4}} \norm{|\nabla u(s)|^2}_{H^{\frac{1}{4}}}$, with $s = \frac{1}{4}$ and $s_1 = s_2 = \frac{15}{16}$. Then $s \leq s_j$ and $s < s_1 + s_2 - \frac{3}{2} = \frac{30}{16} - \frac{24}{16} = \frac{3}{8}$.  We then get
    \begin{align*}
        \norm{|D|^{\alpha + \frac{1}{4}} u(t)}_{L^2} &\lesssim t^{-\frac{\alpha+\frac{1}{4}}{4}} \norm{u_0}_{X_0} + t^{\frac{2 - \alpha}{4}} \norm{u}_{X_T}^3 
        \\ 
        &\hspace{0.3cm} + \int_{0}^t (t-s)^{-\frac{\alpha + 2}{4}} \norm{u(s)}_{H^{1+\frac{15}{16}}}^2 \mathrm{d}s
        \\
        &\lesssim t^{-\frac{\alpha}{4} - \frac{1}{16}}\norm{u_0}_{X_0} + t^{\frac{2-\alpha}{4}}\left[\norm{u}_{L^{\infty}_{(0,T]}H^{\alpha}}^2 + \norm{u}_{X_T}^3\right],
    \end{align*}

    \noindent {\bf Step 2''.}
    Let us apply the operator $|D|^{\alpha + \frac{1}{2}}$, with $\alpha < 2$, to \eqref{mild formulation with time shift}.
    \begin{align*}
        \norm{|D|^{\alpha + \frac{1}{2}} u(t)}_{L^2} \leq &\norm{|D|^{\alpha + \frac{1}{2}}S(t) u_0}_{L^2} +  \int_{0}^t \norm{|D|^{\alpha+\frac{1}{2}}\nabla S(t-s) \cdot \big(|\nabla u(s)|^2 \nabla u(s)\big)}_{L^2} \mathrm{d}s
        \\
        + &\int_{0}^t \norm{|D|^{\alpha} \Delta S(t-s) |D|^{\frac{1}{2}}\big(|\nabla u(s)|^2\big)}_{L^2} \mathrm{d}s.
    \end{align*}
    Using the decay estimates \eqref{decay-est real derivatives}, we have 
    \begin{align*}
        \norm{|D|^{\alpha + \frac{1}{2}} u(t)}_{L^2} &\lesssim t^{-\frac{\alpha+\frac{1}{2}}{4}} \norm{u_0}_{L^2} + \int_{0}^t (t-s)^{-\frac{\alpha + \frac{3}{2}}{4}} \norm{|\nabla u(s)|^2 \nabla u(s)}_{L^{2}} \mathrm{d}s 
        \\ 
        &\hspace{0.3cm} + \int_{0}^t (t-s)^{-\frac{\alpha + 2}{4}} \norm{|\nabla u(s)|^2}_{H^{\frac{1}{2}}} \mathrm{d}s.
    \end{align*}
    Now, let us use Lemma \eqref{prodest} for the terms $I_1 := \int_{0}^t (t-s)^{-\frac{\alpha + \frac{3}{2}}{4}} \norm{|\nabla u(s)|^2 \nabla u(s)}_{L^{2}} \mathrm{d}s$ and $I_2 := \int_{0}^t (t-s)^{-\frac{\alpha + 2}{4}} \norm{|\nabla u(s)|^2}_{H^{\frac{1}{4}}}$.

    For the term $I_1$, we first use the case $s = 0$, for which we take $s_1 = \frac{1}{4} + \varepsilon$ and $s_2 = 1 + \frac{1}{4} - \frac{\varepsilon}{2}$, so that $s_1 + s_2 - \frac{3}{2} = \frac{\varepsilon}{2}$. We get
    \begin{align*}
        I_1 \lesssim \int_{0}^t (t-s)^{-\frac{\alpha + \frac{3}{2}}{4}} \norm{|\nabla u(s)|^2}_{H^{\frac{1}{4} + \varepsilon}} \norm{\nabla u(s)}_{H^{1 + \frac{1}{4} - \frac{\varepsilon}{2}}} \mathrm{d}s.
    \end{align*}
    Then, we use the \ref{Lemma regularity} a second time, in the case $s = \frac{1}{4} + \varepsilon$, for which we take $s_1 = s_2 = 1$, so that $s_1 + s_2 - \frac{3}{2} = \frac{1}{2} > s$ if we take $\varepsilon$ sufficiently small. At the end, we have 
    \begin{align*}
        I_1 \lesssim \int_{0}^t (t-s)^{-\frac{\alpha + \frac{3}{2}}{4}} \norm{\nabla u(s)}_{H^{1 + \frac{1}{4} - \frac{\varepsilon}{2}}}^3 \mathrm{d}s \lesssim t^{\frac{5}{8} - \frac{\alpha}{4}} \norm{u(s)}_{L^{\infty}_{(0,T]}H^{2 + \frac{1}{4}-}}^3
    \end{align*}

    For the term $I_2$, we use the case $s = \frac{1}{2}$, for which we take $s_1 = s_2 = 1 + \frac{\varepsilon}{2}$. Then $s \leq s_j$ and $s < s_1 + s_2 - \frac{3}{2} = \frac{1}{2} + \varepsilon$ and we get
    \begin{align*}
        I_2 \lesssim \int_{0}^t (t-s)^{-\frac{\alpha + 2}{4}} \norm{u(s)}_{H^{2 + \frac{1}{4}-}}^2 \lesssim t^{\frac{2-\alpha}{4}}\norm{u(s)}_{H^{2 + \frac{1}{4}-}}^2.
    \end{align*}
    At the end, we get
    \begin{align*}
        \norm{|D|^{\alpha + \frac{1}{2}} u(t)}_{L^2} \lesssim t^{-\frac{\alpha}{4} - \frac{1}{8}}\norm{u_0}_{X_0} + t^{\frac{2-\alpha}{4}}\norm{u}_{L^{\infty}_{(0,T]}H^{2 + \frac{1}{4}-}}^2 + t^{\frac{5}{8} - \frac{\alpha}{4}}\norm{u}_{L^{\infty}_{(0,T]}H^{2 + \frac{1}{4}-}}^3.
    \end{align*}
    
    Now, we know that $u \in L^{\infty}_{(0,T]} H^{2 + \frac{1}{2}-}(\mathbb{R}^d)$. Then, Step 2'' can be repeated to prove that $u \in L^{\infty}_{(0,T]} H^{2 + \frac{3}{4} -}(\mathbb{R}^d)$.  The rest of the proof is similar to the case $d = 2$, where now $s_0 + 1 = \frac{3}{2} + 1 + \varepsilon < 2 + \frac{3}{4}$ if we take $s_0$ sufficiently small.

\end{remark}

\begin{corollary}\label{corollary regularity}
    Let $u \in X_T$ be a solution of \eqref{dmbe-I} for the initial condition $u_0 \in X_0$. Then $u \in C^{\infty}((0,T];C^{\infty}(R^d))$.
\end{corollary}

\begin{proof}
    From Lemma \ref{Lemma regularity}, we have $u \in L^{\infty}((0,T];C^{\infty}(\mathbb{R}^d))$. Then using the equation \eqref{dmbe-I}, one easily get $u \in C^{0}((0,T];C^{\infty}(\mathbb{R}^d))$. And from equation \eqref{dmbe}, we get $u \in C^{1}((0,T];C^{\infty}(\mathbb{R}^d))$. Then a classical iteration argument, differentiating the equation \eqref{dmbe} with respect to $t$, we obtain $u \in C^{\infty}((0,T];C^{\infty}(\mathbb{R}^d))$.
\end{proof}

\subsection{ Global existence for \eqref{dmbe} }

 From Theorem \ref{thm-wp}\eqref{wpi}, we have the existence of a unique local solution of the MBE equation, $\tilde{u} \in  C^{\infty}((0,T);C^{\infty}(\mathbb{R}^d))$, such that $\tilde{u}(x, 0) = u_0(x)$. In particular, for any $0 < \bar{t} < T$, $\tilde{u}(\bar{t}) \in H^s(\mathbb{R}^d)$, $s\geq 3$. So by \cite[Theorem 3.1]{A15}, there exists a unique global solution $u \in C^0([0, +\infty); H^s(\mathbb{R}^d))$ in the case $d=1,2$. Then from Lemma \ref{Lemma regularity} we get that $u \in C^{\infty}((0,+\infty);C^{\infty}(\mathbb{R}^d))$.

\section{Proof of Theorem \ref{thm-slwp}}\label{lwpstoch}

With the deterministic case complete, we now turn to the stochastic case.  Before we begin, let us fix the notation.  For an $\mathcal F$-measurable function $X : \Omega \rightarrow \R$ and $g: \R \rightarrow \R$, the expectation  $\mathbb E[g(X)]$ is defined as
$$
\mathbb E[g(X)]:= \int_\Omega g(X(\omega)) \, d \mathbb P(\omega) = \int_\R g(x) dF(x)
$$
provided that $\int_\Omega |g(X(\omega))| \, d \mathbb P(\omega) < \infty $, where $F(x) = \mathbb{P} ( X\le x)$ is the distribution function of $X$. 
The $L^p_\omega$-space is then defined via the norm
$$
\norm{X}_{L^p_\omega}=  \left\{\mathbb E\left[ |X|^p \right] \right\}^\frac1p \qquad (1\le p<\infty).
$$
Note that we have the embedding $L^q_\omega \subset L^p_\omega $ whenever $q>p$, i.e., 
$$
\norm{X}_{L^p_\omega} \le \norm{X}_{L^q_\omega}
.$$

\begin{lemma}\label{lm-scov-reg}
Suppose that $T>0$ and $\phi \in HS(L^2; H^s)$ for some $s\in \R$. Then 
\begin{enumerate}[(i)]
\item 
$Z \in C \left( [0,T]; H^s (\R^d)\right) $ almost surely for $\omega \in \Omega$. Moreover, for any  $ p \ge 1$,  there exists $C=C_{p, T} >0$ such that
 $$
\mathbb E \left[   \sup_{0\le t \le T}  \norm{ Z (t)}^p_{H^s}\right] \le C \norm{\phi}^p_{HS(L^2: H^s)}.
$$
 
 \item $Z \in L^q_T\dot W^{1,r} (\R^d)$ for $q\ge 1$ and  $r\ge 2$, almost surely for $\omega \in \Omega$. Moreover,  for any  $ p \ge 1$, there exists $C=C_{q, p, r, T} >0$ such that
$$
\mathbb E \left[     \norm{  Z }^p_{ L^q_T\dot W^{1,r}}\right] \le C  \norm{\phi}^p_{HS(L^2: L^2)}.
$$
\end{enumerate}

\end{lemma}
\begin{proof}
Part (i) is contained in Theorem 6.10 of \cite{PZ14}.  So we only prove (ii). 

From \eqref{stconv}, we see that $ \nabla  Z(x,t)$ is a linear combination of independent Wiener integrals. It is a mean zero Gaussian random process with variance 
\begin{equation}\label{var}
\begin{split}
\sigma(x,t) &=\norm{  \nabla  Z(x,t)}^2_{L^{2}_\omega} =\mathbb E \left[   | \nabla  Z(x,t) |^2 \right]
\\
&=\sum_{k\in \N} \int_0^t  |S(t-s) \nabla  (\phi e_k)(x)   |^2\;  ds,
\end{split}
\end{equation}
where we used the It\^{o} isometry to obtain the last equality.

Now, recall that for mean zero Gaussian random $X$ variable with variance $\sigma$, we have
\begin{equation}\label{varf}
\mathbb E \left[   |X|^{2j} \right]=j! \sigma^j.
\end{equation}
So for a  positive integer $p$, we have by the embedding $L^{2p}_\omega \subset L^p_\omega $, \eqref{varf} and \eqref{var},
\begin{equation}
\label{ExpZ}
\begin{split}
\norm{ \nabla Z(x,t)}_{L^p_\omega}  &\le  \norm{  \nabla  Z(x,t)}_{L^{2p}_\omega} 
\\
&\le  C_p \norm{  \nabla  Z(x,t)}_{L^{2}_\omega} 
\\
&\sim  \left[\sum_{k\in \N} \int_0^t  |S(t-s)  \nabla   (\phi e_k)(x)   |^2\;  ds\right]^\frac12
\\
&= \norm{ S(t-\cdot)  \nabla   (\phi e_k)(x) }_{l^2_k L^2_{[0,t]}}.  
\end{split}
\end{equation}
%%%%%%%%%%%%%%%%%%%%%%%%%%%%%%%%%%%%

%%%%%%%%%%%%%%%%%%%%%%%%%%%%%%%%%

Assume $p \ge \max( q, r)$. By \eqref{ExpZ} and the decay estimate \eqref{decay-est real derivatives}
\begin{align*}
 \norm{ \nabla  Z( t)}_{ L^p_\omega L^q_T L_x^r}  &\le   \norm{ \norm{  \nabla  Z (x,t)}_{L^p_\omega}}_{  L^q_T L_x^r}
\\
&\le C_p \norm{ \norm{  S(t-\cdot)  \nabla   (\phi e_k)(x) }_{l^2_k L^2_{[0,t]}} }_{ L^q_T L_x^r}
\\
&\le C_{p} \norm{ \norm{  \norm{S(t-\cdot)  \nabla   (\phi e_k) }_{L_x^r} }_{l^2_k L^2_{[0,t]}} }_{L^q_T }
\\
&\le C_{p, q,r} \norm{  \norm{  (t-\cdot)^{ \frac d{4r}-\frac {2+d}8} \norm{(\phi e_k)}_{L_x^2} }_{l^2_k L^2_{[0,t]}} }_{L^q_T }
\\
&= C_{p, q,r} \norm{   t^{ \frac d{4r}-\frac d8 + \frac 14} \norm{  \norm{ \phi e_k}_{L_x^2} }_{l^2_k } }_{L^q_T }
\\
&\le  C_{p,q,r}  T^{\frac 1q+ \frac d{4r}-\frac d8 + \frac 14} \norm{ \phi }_{ HS(L^2: L^2 )}.
\end{align*}

\end{proof}

\subsection{ Local existence for \eqref{ismbe} }
For simplicity, we take $\alpha_1 = \alpha_3 = 1$, $\alpha_2 = \alpha_4 = 0$. However, the argument works for any $\alpha_j \ge 0$ ($j=1, 2, 3, 4$). So $N(u)$ becomes
$$N(u)=\nabla\cdot\left(|\nabla u|^2\nabla u\right) - \Delta u.$$

To prove local existence of solution \eqref{smbe-I-resd} we apply a fixed point argument in the ball
 $
B_{ R, T}= \left\{  v \in X_T : \  \norm{v}_{X_T}  \le 2 R \right\},
$
where
$$X_T=  C \left( [0,T]; X_0 \right).$$

For any $v \in B_{ R, T}$, 
we define the map $\Gamma= \Gamma_{u_0, Z}$,
by
$$
\Gamma v(t)=   S(t)  u_0
	+  \int_0^t S(t-s) N(v+Z)  (s)\; ds,
$$
where
$$
N(v+Z) = \nabla \cdot \left ( |\nabla (v+ Z) |^2  \nabla (v+Z) \right)  -\Delta (v+Z).
$$

First note that by Lemma \ref{lm-scov-reg}, 
$$  Z \in C \left( [0,T]; L^2 \right)  \cap   L^q_T\dot W^{1,4} (\R^d) \qquad (q\ge 1)$$ almost surely for $\omega \in \Omega$. In particular, there exists for  $A=A(q, T)>0$ such that
$$
 \norm{  Z }_{ L^\infty_T L^2_x}+ \norm{  Z }_{ L^q_T\dot W^{1,4}} \le A $$
almost surely for $\omega \in \Omega$.  Thus, local existence of a unique solution in $B_{ R, T}$ follows if we prove 
\begin{align}
\label{Gmap}
\norm{\Gamma v}_{X_T} & \le   2R 
\\
\label{Gmap-d}
\norm{\Gamma v- \Gamma w}_{X_T}  < &\norm{v- w}_{X_T}
\end{align}
 for $v, w\in  B_{ R, T}$ and sufficiently small $T = T(R,A)$.

\subsection*{Proof of \eqref{Gmap}}
Fix $\rho>8$. Following the same argument as in the deterministic case, we obtain
\begin{equation*}
\begin{split}
\norm{\Gamma v(t)}_{L^2_x} \le & \norm{S(t) u_0}_{ L^2_x} \\
& +  \int_0^t \left\{ (t-s)^{-\frac 12 }  \norm{(v+Z)(s)}_{L_x^2}  +   (t-s)^{-\frac {4+d}{16} }  \norm{  \nabla (v+ Z )(s)}^3_{L_x^{4}}  \right\}
\; ds 
\\ 
%\le & \norm{ u_0}_{ L^2_x}+ C \left[T^\frac12 (\norm{v}_{L^\infty_TL_x^2} +\norm{Z}_{L^\infty_TL_x^2} ) +  \int_0^t  (t-s)^{-\frac {4+d}{16} }  \norm{ \nabla v (s)}^3_{L_x^4} \; ds +  \int_0^t  (t-s)^{-\frac {4+d}{16} }  \norm{ \nabla Z (s)}^3_{L_x^4} \; ds \right]
%\\ 
\le & \norm{ u_0}_{ L^2_x}+C \left[ T^\frac12 (\norm{v}_{L^\infty_TL_x^2} +\norm{Z}_{L^\infty_TL_x^2} )\right] \\
& + C\left[ T^{\frac {12-d}{16} }  \norm{ \nabla v }^3_{L_T^\infty L_x^4} + T^{\frac {12-d}{16} -\frac1\rho}  \norm{ \nabla Z }^3_{L_T^{3\rho} L^4_x} \right].
\end{split}
\end{equation*}
In the last inequality, we used H\"{o}lder inequality in time to obtain the last term.  Similarly, we may also obtain
\begin{align*}
\norm{ \nabla \Gamma v(t)}_{L^4_x} \le & \norm{   \nabla S(t) u_0}_{ L^4_x} \\
& + \int_0^t \left\{(t-s)^{-\frac {12+d}{16} }  \norm{(v+Z)(s)}_{L_x^2}  + (t-s)^{-\frac {4+d}{8}}  \norm{  \nabla (v+ Z )(s)}^3_{L_x^{4}}   \right\}
\; ds 
\\ 
\le & \norm{   \nabla u_0}_{ L^4_x}  + C \left[T^{\frac {4-d}{16} }  (\norm{v}_{L^\infty_TL_x^2} +\norm{Z}_{L^\infty_TL_x^2} ) \right] \\
& + C\left[ T^{\frac {4-d}{8} }   \norm{v_x}^3_{L^\infty_TL_x^4}  + T^{ \frac {4-d}{8}-\frac1\rho} \norm{ \nabla Z }^3_{L_T^{3\rho}  L^{4}_x} \right].
\end{align*}

From the above estimates, we conclude that for some $\theta>0$,
 \begin{align*}
\norm{\Gamma v}_{X_T} &\le \norm{u_0}_{ X_0} + C  T^\theta \left[   \norm{v}_{ X_T}+  \norm{v}^3_{ X_T} + \norm{Z}_{ L^\infty_T L^2_x} +\norm{Z}^3_{ L^{3\rho}_T \dot W^{1,4}}  \right]
\\
&\le  R + C T^\theta  \left(2R+ 8 R^3 + A+A^3\right)
\\
&\le  2R ,
 \end{align*}
provided that
$$
   T\le \left[\frac{R}{ C\left(2R+ 8R^3 + A+A^3\right)} \right]^\frac1\theta  .$$
This proves \eqref{Gmap}.
%\vspace{2mm}

\subsection*{Proof of \eqref{Gmap-d}}
Begin by writing
\begin{align*}
 N(v+Z) -N(w+Z) 
&=  \nabla \cdot \left[|\nabla v|^2\nabla ( v -w) + \left(\nabla ( v -w) \cdot \nabla ( v +w) \right) \nabla w \right] 
\\
& \quad +  \nabla \cdot  \left[ |\nabla Z|^2\nabla( v-w) +  \left(\nabla (v-w)\cdot \nabla (v+w)\right) \nabla Z\right]
\\
 & \quad + 2 \nabla \cdot \left[   \left(\nabla (v-w)\cdot   \nabla Z \right) \nabla v + \left(\nabla w \cdot   \nabla Z \right) \nabla (v-w)\right] \\ & \quad + 2 \nabla \cdot \left[ \left( \nabla (v-w)\cdot \nabla Z \right)  \nabla Z \right] -\Delta (v-w).
\end{align*}
Using this identity, we see that
\begin{align*}
\norm{\Gamma v(t)-\Gamma w(t) }_{L^2_x} \le & C \int_0^t (t-s)^{-\frac 12 }  \norm{(v-w)(s)}_{L_x^2} \; ds \\
& + C \int_0^t \left\{(t-s)^{-\frac {4+d}{16} }  \left[  \norm{ \nabla  v(s)}^2_{L_x^{4}} + \norm{ \nabla w(s)}^2_{L_x^{4}}+ \norm{ \nabla  Z(s)}^2_{L_x^{4}}\right]\right. \times \\
& \quad \times \left.\norm{ \nabla  (v- w )(s)}_{L_x^{4}}  \right\}
\; ds 
\\ 
&   \le C T^\frac12 \norm{v-w}_{L^\infty_TL_x^2} \\
& +  C \left[ T^{\frac {12-d}{16} }  \left( \norm{ \nabla v }^2_{L_T^\infty L_x^4} + \norm{ \nabla w }^2_{L_T^\infty L_x^4} \right)+  T^{\frac {12-d}{16} -\frac1\rho} \norm{ \nabla Z }^2_{L_T^{2\rho } L_x^{4}} \right] \times \\
& \quad \times \norm{ \nabla (v-w) }_{L_T^\infty L_x^4}.
\end{align*}
A similar computation yields 
\begin{align*}
\norm{\nabla \left( \Gamma v(t)-\Gamma w(t) \right) }_{L^4_x} \le & C T^{\frac {4-d}{16} } \norm{v-w}_{L^\infty_TL_x^2} \\
& + C \left[ T^{\frac {4-d}{8} } \left( \norm{ \nabla v }^2_{L_T^\infty L_x^4} + \norm{ \nabla w }^2_{L_T^\infty L_x^4} \right) \right. \\
& + \left. T^{ \frac {4-d}{8}-\frac1\rho}  \norm{ \nabla Z }^2_{L_T^{2 \rho} L_x^{4}} \right] \norm{ \nabla (v-w) }_{L_T^\infty L_x^4}.
\end{align*}
Thus, 
\begin{align*}
\norm{\Gamma v-\Gamma w}_{X_T} &\le C  T^\theta( R^2 + A^2 ) \norm{ v- w}_{X_T} 
<  \norm{ v- w}_{X_T} 
 \end{align*}
provided that
$$C T^\theta ( R^2 + A^2 ) < \left(\frac{1}{C \left(R^2 + A^2\right)}\right)^{\frac{1}{\theta}}. $$
This completes the proof of \eqref{Gmap-d}.

\section{Proof of Lemma \ref{lm-ap}}
Suppose that $v$ is a smooth solution of \eqref{smbe-I-resd} on $[0, T]$.  Then $v$ satisfies
\begin{equation}\label{smbe-I-resd1}
\left\{
\begin{aligned}
&v_{t} +\Delta^2  v =  N(v+Z)  \\
& v(x,0) = u_0(x),
\end{aligned}
\right.
\end{equation}
where (by assumption $\alpha_1 = \alpha_3 = 1$, $\alpha_2 = \alpha_4 = 0$)
$$
N(v+Z) = \nabla \cdot \left ( |\nabla (v+ Z) |^2  \nabla (v+Z) \right)  -\Delta (v+Z).
$$

Let
\begin{align*}
M(v(t))=  \frac12 \norm{v(t)}_{L_x^2}^2 +  \int_0^t \left[ \norm{ \Delta v (s)}_{L_x^2}^2 +  \norm{ \nabla v (s)}_{L_x^4}^4\right ] \, ds.
\end{align*}
Then by \eqref{smbe-I-resd1} and Plancherel,
\begin{align*}
\frac {dM (v(t))}{dt} &=\int_{\R^d} v v_t +    (\Delta v)^2+ |\nabla v|^4 \, dx 
\\
= &  \int_{\R^d } v \left[ v_t +  \Delta^2 v - \nabla \cdot ( |\nabla v|^2 \nabla v) \right] \, dx
\\
= &\int_{\R^d}  v\nabla \cdot \left[ 2 (\nabla v \cdot \nabla Z) ( \nabla v + \nabla Z)  + |\nabla v|^2 \nabla Z +  |\nabla Z|^2 \nabla v+  |\nabla Z|^2 \nabla Z \right] \, dx \\
& -\int_{\R^d} v \Delta  (v+Z) \, dx 
\\
= & \int_{\R^d}  \nabla v \cdot \left[- 2 (\nabla v \cdot \nabla Z) ( \nabla v + \nabla Z)  - |\nabla v|^2 \nabla Z -  |\nabla Z|^2 \nabla v-  |\nabla Z|^2 \nabla Z \right] \, dx \\
& + \int_{\R^d} [- (v+Z) ] \Delta v \, dx .
\end{align*}
 By Young's inequality,  the integrad in the first integral on the right is bounded by 
 \begin{align*}
   |\nabla v|^4   + c  |\nabla Z|^4   
 \end{align*}
 for some $c\gg 1$. For instance, 
\begin{align*}
\vert \nabla v \cdot \left[  (\nabla v \cdot \nabla Z) ( \nabla v + \nabla Z) \right) \vert \le & \frac 34 \left( \frac12 |\nabla v|\right)^4 +  \frac 14 \left(8|\nabla Z| \right)^4 \\
& + \frac 12 \left( \frac12 |\nabla v|\right)^4 + \frac 12 \left( 2|\nabla Z| \right)^4.
\end{align*}
 Therefore, the first integral on the right 
  is bounded by 
  \begin{align*}
  \int_{\R} \left(   |\nabla v|^4   + c  |\nabla Z|^4    \right)\, dx 
 \end{align*}
 Similarly, the second integral is bounded by 
 \begin{align*}
  \int_{\R} \left( \frac 12v^2   +  (\Delta v)^2+ \frac12  Z^2 \right)\, dx 
 \end{align*}
Thus, 
\begin{align*}
  \frac {dM (v(t))}{dt} & \le  M (v(t))+  c \norm{ \nabla Z(t)}^4_{L_x^4} +   \frac12 \norm{Z(t)}^2_{L_x^2}.
 \end{align*}
Applying Gr\"onwall's inequality, it is evident that
\begin{equation}\label{gronwall}\begin{aligned}
M( v(t))  &\le  M( v(0)) e^{t} + C_T(e^t-1)
\\
&\le \left[ \norm{u_0}^2_{L_x^2}  + C_T\right] e^{t}  ,
 \end{aligned}\end{equation}
where 
$$C_T= \sup_{0\le t\le T} \left( c \norm{ \nabla Z(t)}^4_{L_x^4} +   \frac12 \norm{Z(t)}^2_{L_x^2}\right).$$
This proves Lemma \ref{lm-ap}.

% \vspace{2mm}
\section{On Coarsening in MBE Models}\label{coarsening}

%One of the major technical problems in MBE is the phenomenon of coarsening.  Ideally, surfaces generated by MBE should be perfectly flat.  However, because mounds and other imperfections have been observed experimentally, it is known that the manufactured surfaces cannot be expected to be perfectly flat.  Because of this, there is a great interest in estimating the size of the imperfections, as this could affect the performance of the manufactured component.  As is often the case in engineering, the size of the deviations are typically measured using a root-mean-square calculation, which is conveniently captured by the $L^{2}$ norm of solutions.  This leads us to the following definition:
Continuing the discussion of coarseness in the introduction, let us state the formal definition of the coarsness of a solution.
\begin{definition}
    Let $u(x,t)$ be a solution to any partial differential equation modeling MBE.  We define the coarseness $C(t)$ of $u$ at time $t$ by
    \[
    C(t) = \| u(\cdot, t) \|_{L^{2}_{x}}.
    \]
\end{definition}
One of the most well-known asymptotic estimates for the growth rate of the coarseness is due to Rost and Krug, who showed that solutions to their model satisfy the estimate
\[
C(t) \leq Ct^{\frac{1}{2}}
\]
in the deterministic case.  Interestingly, this estimate was obtained using only basic calculus and some mild assumptions on solutions based on physical principles.  With this fact in mind, we observe that Theorem 3.1 of \cite{A15} implies that the solutions obtained in Theorem \ref{thm-wp} must obey the a priori estimate
\[
C(t) \leq C(0) e^{Kt},
\]
where $K$ is a non-negative constant which depends on the coefficients $\alpha_{1},\ldots, \alpha_{4}$.  On the surface, this would appear to be a much weaker result.  However, it should be noted that in a typical situation encountered in manufacturing, we start with a flat surface, so that $u(0) = 0$.  This suggests that the solutions constructed in Theorem \ref{thm-wp} should exhibit no coarsening.  Since the model of Rost and Krug is a special case of equation \eqref{mbe}, then the result should hold for that model, as well.  Moreover, the uniqueness of solutions in the space $X_{T}$ implies that any solution corresponding to zero initial data which exhibits coarsening cannot be in $X_{T}$.  The only way that this can happen is if either
\[
\| u(t_{0}) \|_{L^{2}_{x}} = + \infty \quad \text{or} \quad \| \nabla u(t_{0}) \|_{L^{4}_{x}} = + \infty
\]
for some $t_{0} > 0$.  The former case is particularly troubling, as such solutions would either have to go to infinity in the continuous case, or be discontinuous.  Since both cases are physically unrealistic, this suggests that deterministic models should not lead to solutions which exhibit coarsening, which goes against the conventional wisdom and experimental results in the study of MBE.

To reconcile this discrepancy, we remark that in physically realistic situations, signals are often noisy.  Thus, an accurate model of epitaxial dynamics should include noise.  Based on the results of Theorem \ref{thm-slwp} and Lemma \ref{lm-ap} (in particular, equation \eqref{gronwall}), global solutions perturbed by a noise should satisfy the coarsening estimate
\[
C^{2}(t) \le C^{2}(0)e^{t} + C_{T} (e^{t} - 1).
\]
For initially flat surfaces, this reduces to
\[
C^{2}(t) \le C_{T} (e^{t} - 1).
\]
As is well known,
\[
e^{t} - 1 = \mathcal{O}(t)
\]
as $t \rightarrow 0$.  We may thus conclude that
\[
C(t) \le Ct^{\frac{1}{2}},
\]
for short times, so that we recover Rost and Krug's result.  This suggests that any model of MBE must contain a noise term to exhibit coarsening.  Subsequently, it also suggests that by improving manufacturing methods, it may be possible to reduce the coarsening of surfaces.

%%%%%%%%%%%%%%%%%%%%%%%%%%%%%%%%%%%%%%%
%\vspace{5mm}

\section*{Acknowledgments}
L.E. and A. T. are supported by the Faculty Development Competitive Research Grants Program 2022-2024, Nazarbayev University: {Nonlinear Partial Differential Equations in Material Science  (Ref. 11022021FD2929)}

%\vspace{5mm}

%%%%%%%%%%%%%%%%%%%%%%%%%%%%%%%%%%%%%%%%%%%%%%%%%%%%%%

\end{document}